\documentclass{amsart}[10pt]

\usepackage{amsthm}
\usepackage{graphicx}
\usepackage{amssymb}
 \calclayout

\newtheorem{teo}{Theorem}[section]

\theoremstyle{definition}
\newtheorem{defi}[teo]{Definition}
\newtheorem{rem}[teo]{Remark}

\theoremstyle{remark}
\newtheorem{prof}[teo]{Proof of}

\def\mc{\mathbb{C}}
\def\mr{\mathbb{R}}

\begin{document}
\title{A calculus for branched spines of 3-manifolds}

\author[Costantino]{Francesco Costantino}

\address{Scuola Normale Superiore\\
  Piazza dei Cavalieri 7\\
  56126 Pisa, Italy}
\email{f.costantino@sns.it}


\begin{abstract}
We establish a calculus for branched spines of $3$-manifolds by
means of branched Matveev-Piergallini moves and branched
bubble-moves. We briefly indicate some of its possible
applications in the study and definition of State-Sum Quantum
Invariants.
\end{abstract}

\maketitle

\tableofcontents
\section{Introduction}

Since the establishment of Matveev-Piergallini calculus, simple
spines of $3$-manifolds have been one of the most powerful tools
to study these spaces. They allowed a combinatorial approach to
many fundamental topics as the study of State-Sum Quantum
Invariants and the study of complexity of $3$-manifolds.

Branched spines of $3$-manifolds, which could be viewed as the
smoothed version of simple spines, were first introduced and
studied by Benedetti and Petronio in \cite{BP}. Among other
substantial results, in this book, Benedetti and Petronio showed
that each $3$-manifold has a branched spines and identified the
topological structure encoded by these objects on the ambient
manifolds as a particular class of non-singular vector fields they
called {\it Concave Transversing Fields}. At a more rough level,
it can be showed that branched spines can be used to represent the
$Spin^c$-structures on the ambient manifolds.

In the present work, instead of viewing branched spines as a tool
to represent $3$-manifolds equipped with additional structures, we
will use these objects to re-obtain a representation theory of
naked $3$-manifolds. To clarify the reason why we are interested
in such a representation theory we notice that a branched spine is
much less symmetric than a non-branched one. Moreover, a branched
spine is dual to a triangulation of the ambient manifold whose
``abstract'' tetrahedra can be canonically ``parametrized'' by the
standard simplex $\Delta(v_0,v_1,v_2,v_3)$ in $\mr^3$: indeed,
using the branching, we can canonically identify in each dual
tetrahedron the vertex corresponding to each $v_i$, $i=0,\ldots
3$.

This kind of branched triangulations underlie the
definition of the Quantum Invariants obtained as State-Sum 
and this is the main motivation for the present note. Indeed, 
using these objects, S.Baseilhac and R.Benedetti constructed in 
\cite{BB2} and \cite{BB3} (see also \cite{BB}) the so-called classical and 
quantum dilogarithmic invariants for $3$-manifolds equipped with principal flat 
$PSL(2,\mathbb{C})$-bundles. 

%
%
In the present paper we show, in particular, that two branchings
on a given triangulation can be connected by means of a finite
sequence of basic ($2\to 3$, and ``bubble'') branched modification
of branched triangulations. With the terminology of
Baseilhac-Benedetti, this implies that {\it the branching state
sum invariance is a consequence of its $\it transit$ invariance}.


Before plunging into the details, let us note that what we
call ``calculus for branched spines'' should be better called
``calculus for branched skeleta'' in the sense that we show that
two branched spines of the same $3$-manifold are connected by a
sequence of moves which contains an (algebraic) number of
bubble-moves equal to zero; hence, during the sequence, the spines
could be transformed into spines of the manifold obtained by
puncturing the ambient manifold, which are also called skeleta of
the initial manifold. Fortunately this causes no harm since the
initial ambient manifold can always be canonically reconstructed
from a skeleton. In terms of dual triangulations this means the
natural fact that the number of vertices is not fixed.

{\bf Acknowledgements.} The author wishes to warmly thank Riccardo
Benedetti and Stephane Baseilhac for their encouraging comments
and illuminating critics.

\section{Preliminaries}\label{sub:prelim}

In this section we recall the notion of branched polyhedron and
some basic facts about branched spines of $3$-manifolds. From now
on, we will deal only with polyhedra which have the property of
containing only regions which are orientable surfaces (this is due
to our definition of branching, see Definition
\ref{branchingcondition}) and with oriented $3$-manifolds.

\begin{figure}
  \centerline{\includegraphics[width=8.4cm]{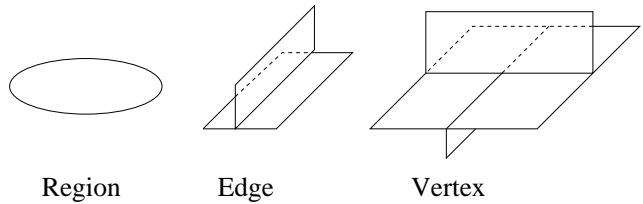}}
  \caption{The local models of a simple polyhedron.  }
  \label{fig:singularityinspine}
\end{figure}

A {\it simple polyhedron} is a finite polyhedron of dimension $2$
whose local model are the three shown in Figure
\ref{fig:singularityinspine}. An {\it embedded spine} of a
$3$-manifold $M$ is a simple polyhedron $P$ embedded in $M$ in a
locally flat way (i.e. so that there exist local charts as those
shown in Figure \ref{fig:singularityinspine}) so that $M$ retracts
on $P$ (if $\partial M\neq \emptyset$) or $M- \{one\  point \}$ retracts 
on $P$ (if $\partial M=\emptyset$).

The {\it singular set} of a simple polyhedron $P$, denoted by
$Sing(P)$, is the union of the edges and vertices (see Figure
\ref{fig:singularityinspine}). The {\it regions} of a simple
polyhedron $P$ are the connected components of the complement of a
small open regular neighborhood of $Sing(P)$.

\begin{figure} [h!]
   \centerline{\includegraphics[width=9.4cm]{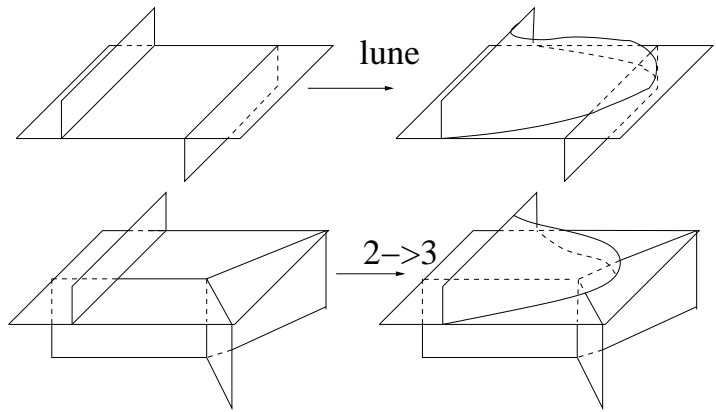}}

  \caption{In this figure we show the basic moves for polyhedra. Note
    that both moves create a new region: the small disc entirely
    contained in the left part of the figures. }\label{2mosse}
\end{figure}

 A simple polyhedron is said to be {\it standard} if its regions are all discs, its
singular set is connected and contains at least a vertex.
 It is well known that any
$3$-manifold admits a spine, even a standard one, moreover the
following holds:
\begin{teo}[Matveev-Piergallini]
Any two standard spines of the same $3$-manifold are connected by
means of a suitable sequence of local moves (and their inverses)
as the one shown in the lower part of Figure \ref{2mosse} and
called the $2\to 3$-move. More in general, two simple spines of
the same $3$-manifold are connected by a suitable sequence of this
move and moves of the type shown in the upper part of Figure
\ref{2mosse} and their inverses; this last move is called the
\rm{lune-move} or $0\to 2$-move.
\end{teo}
\begin{figure} [h!]
   \centerline{\includegraphics[width=13cm]{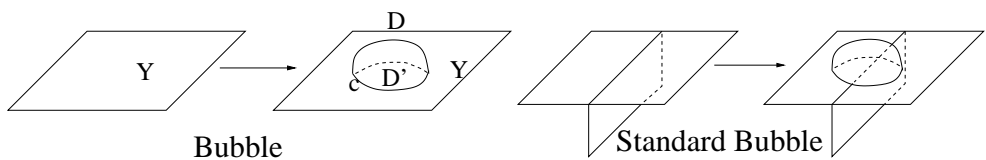}}

  \caption{In the left part of the figure we show the bubble-move; it
    can be interpreted as the gluing of a disc $D$ along its boundary 
to a simple closed curve $c$ contained in the interior of a region $Y$
and bounding a small disc $D'$.
    The result of the application of such a move to a standard polyhedron is
    not a standard polyhedron, so often one uses its standard version, shown in tyhe right part of the picture.
    }\label{fig:bubble}
\end{figure}

The above Theorem is the well known Calculus for spines of
$3$-manifolds proved independently by Matveev and Piergallini
respectively in \cite{Ma} and \cite{Pi}.

Another local move we will use is called {\it bubble-move} it is
applied in the interior of a region
and its effect is drawn in Figure \ref{fig:bubble}.

 Given a simple polyhedron $P$ we define the notion of
{\it branching} on it as follows:
\begin{defi}[Branching condition]\label{branchingcondition}
A branching $b$ on $P$ is a choice of an orientation of each
region of $P$ such that no edge of the singular set of $P$ is
induced three times the same orientation by the regions containing
it. 
\end{defi}

\begin{rem}
This definition corresponds to the definition of ``orientable
branching'' given in \cite{BP}.
\end{rem}

\begin{figure} [h!]
   \centerline{\includegraphics[width=11.4cm]{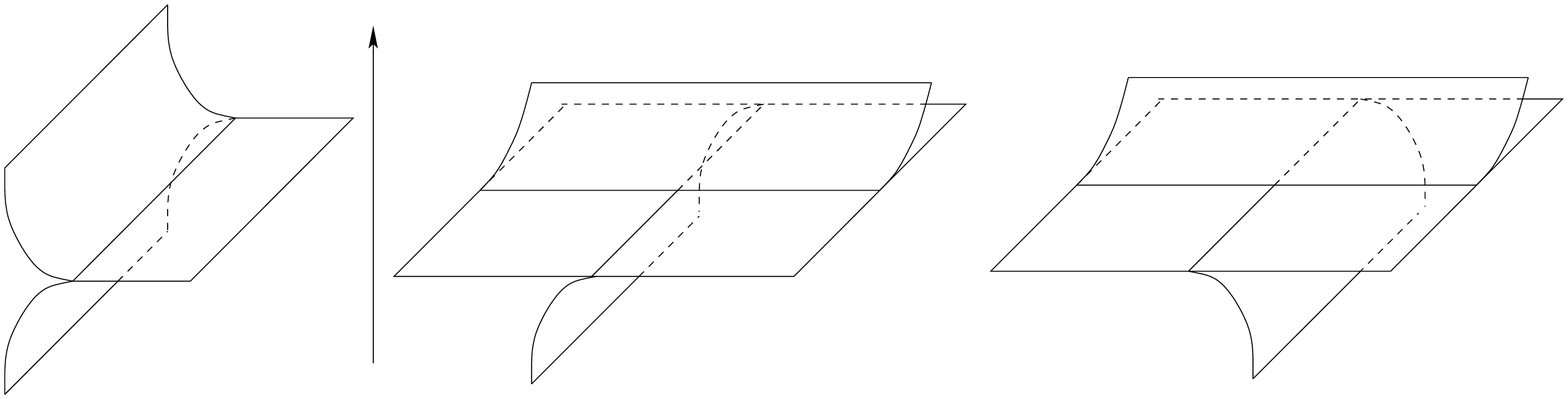}}
  \caption{How a branching allows a smoothing of the polyhedron: the
  regions are oriented using the right-hand rule, the upward direction and
  the orientation of the ambient $3$-manifold (coinciding with the
  standard one of the chart depicted here).}\label{branching}
\end{figure}

Not all the simple polyhedra admit a branching and, on the
contrary, there are some which admit more than one (see \cite{BP},
Chapter III). We
will say that a polyhedron is {\it branchable} if it admits a
branching and we will call {\it branched polyhedron} a pair
$(P,b)$ where $b$ is a branching on the polyhedron $P$.
\begin{defi}
Let $P$ be a spine of an oriented $3$-manifold $M$; $P$ is said to
be {\it branchable} if the underlying polyhedron is. We call {\it
branched spine} of $M$ the pair $(P,b)$ where $P$ is a spine and
$b$ is a branching on the underlying polyhedron. When this will
not cause any confusion, we will not specify the branching $b$ and
we will simply write $P$.
\end{defi}
A branching on a simple polyhedron allows us to smoothen its
singularities and equip it with a smooth structure as shown in
Figure \ref{branching}.

If $P$ is a branched spine of a $3$-manifold, and we apply to it a basic
move (one of the moves of Figure \ref{2mosse}), we get another
spine $P'$ of the same manifold containing one region more than $P$. 
Moreover, each region
of $P$ naturally corresponds to a region of $P'$ and the
region of $P'$ which does not correspond to one of $P$ is the small 
disc created by the move (see Figure
\ref{2mosse}). Hence the branching on $P$ induces a choice of an orientation
on each region
of $P'$ except on that disc and these orientations satisfy the branching condition
on all the edges of $P'$ not touched by that disc. Analogously, if $P'$ is
obtained from $P$ through the inverse of a basic move, then each region of
$P'$ corresponds to a region of $P$ and hence the branching on $P$
induces an orientation on each region of $P'$.

\begin{defi}\label{branchedmove}
A basic move $P\rightarrow P'$ applied on a branched polyhedron $P$ is called {\it
  branchable} if
it is possible to choose an orientation on the disc created by the move which,
together with the orientations on the regions of $P'$ induced by the branching of $P$,
defines a branching on $P'$. Analogously, the inverse of a basic move
  applied to $P$ is branchable if the orientations induced by the
  branching of $P$ on the regions of $P'$ define a branching.
\end{defi}

\begin{figure} [h!]
  \centerline{\includegraphics[width=0.1cm]{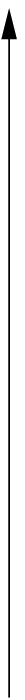}\includegraphics[width=8.4cm]{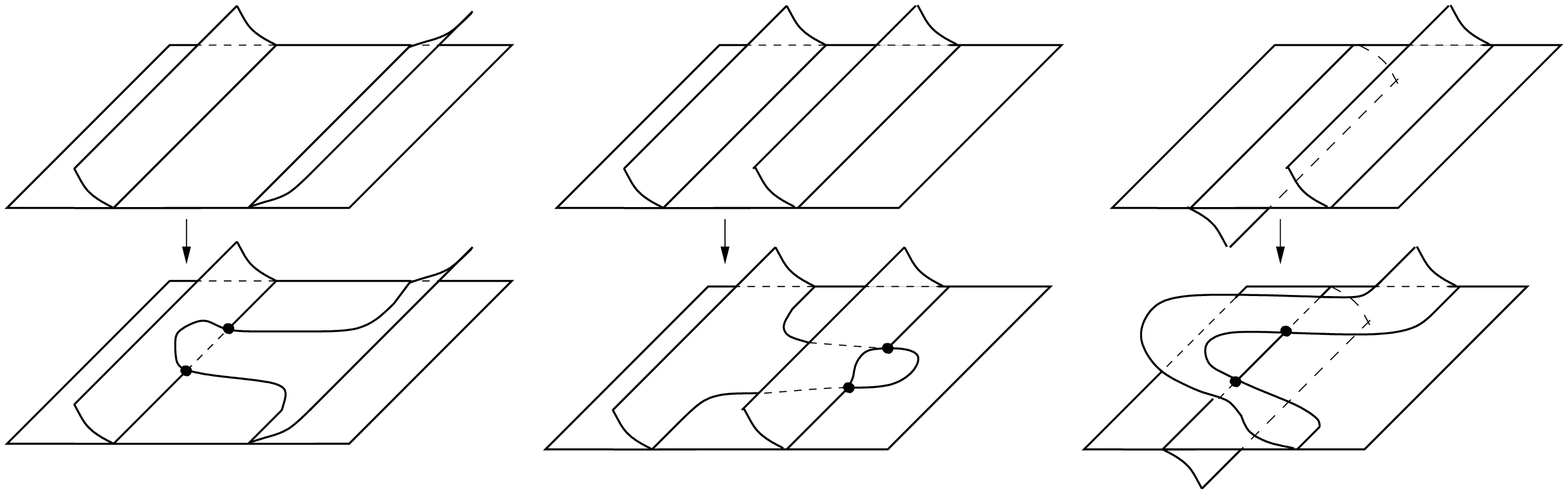}}
 \centerline{\includegraphics[width=0.1cm]{verticaldir.eps}\includegraphics[width=8.4cm]{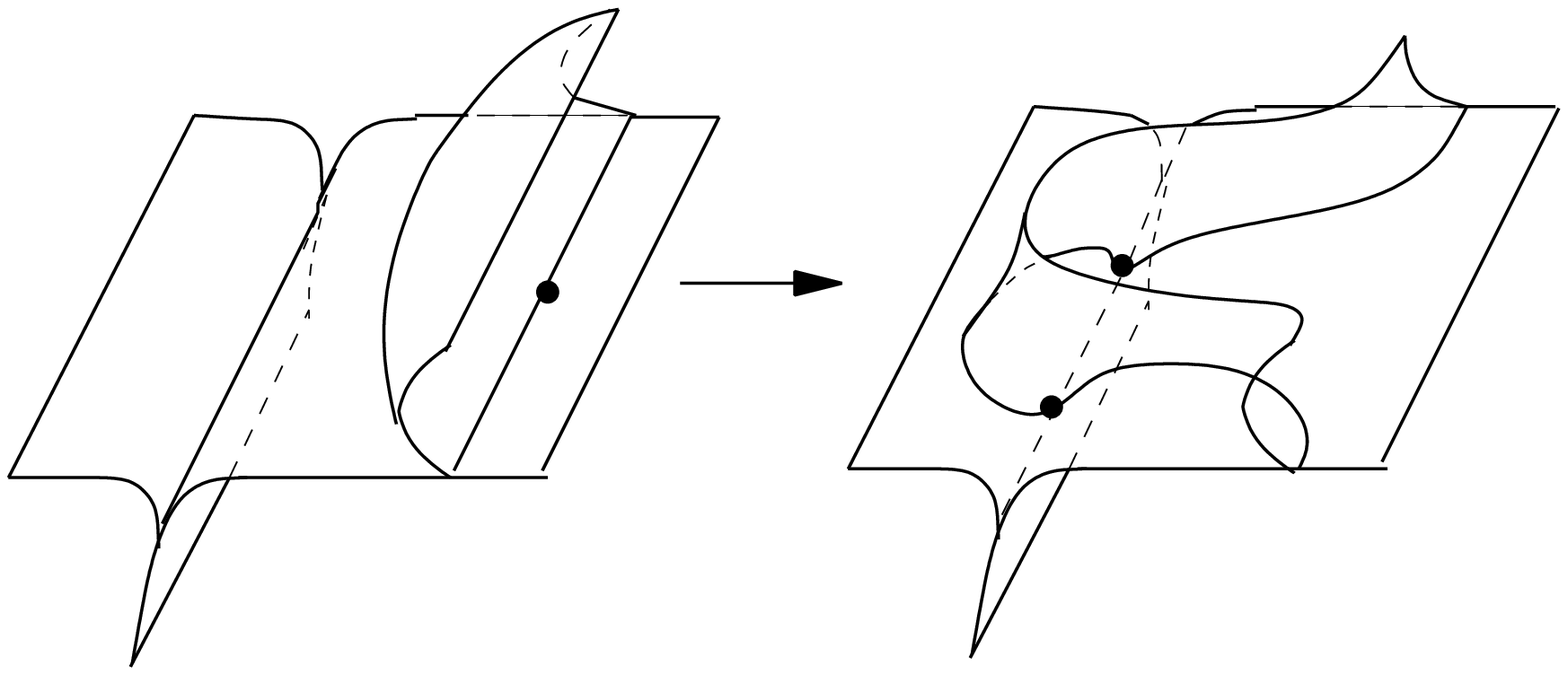}}

  \caption{In the upper part of this figure we show the three branched
  versions of the lune-move called
  ``sliding''-moves. In the bottom part we show the version
  called ``bumping''-move. The arrow on the left indicates the
  vertical direction.}\label{branched lune}
\end{figure}

\begin{figure} [h!]
\centerline{\includegraphics[width=0.1cm]{verticaldir.eps}\includegraphics[width=11.4cm]{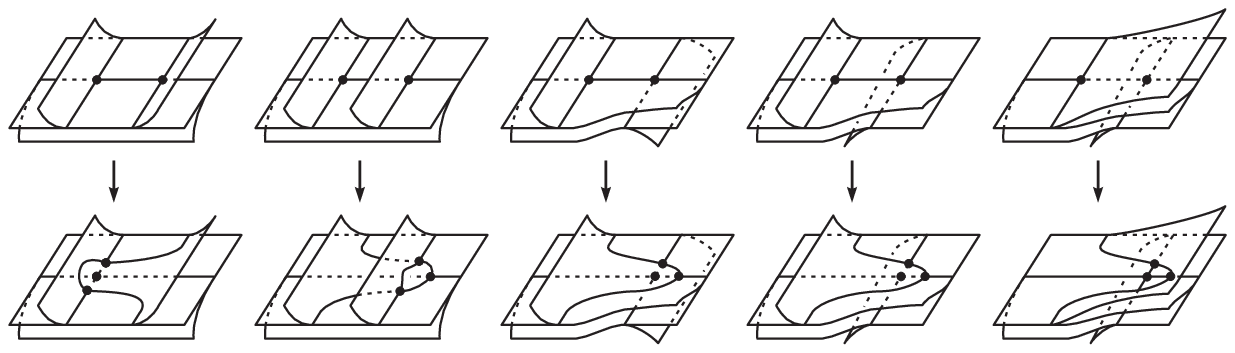}}
\centerline{\includegraphics[width=0.1cm]{verticaldir.eps}\includegraphics[width=8.4cm]{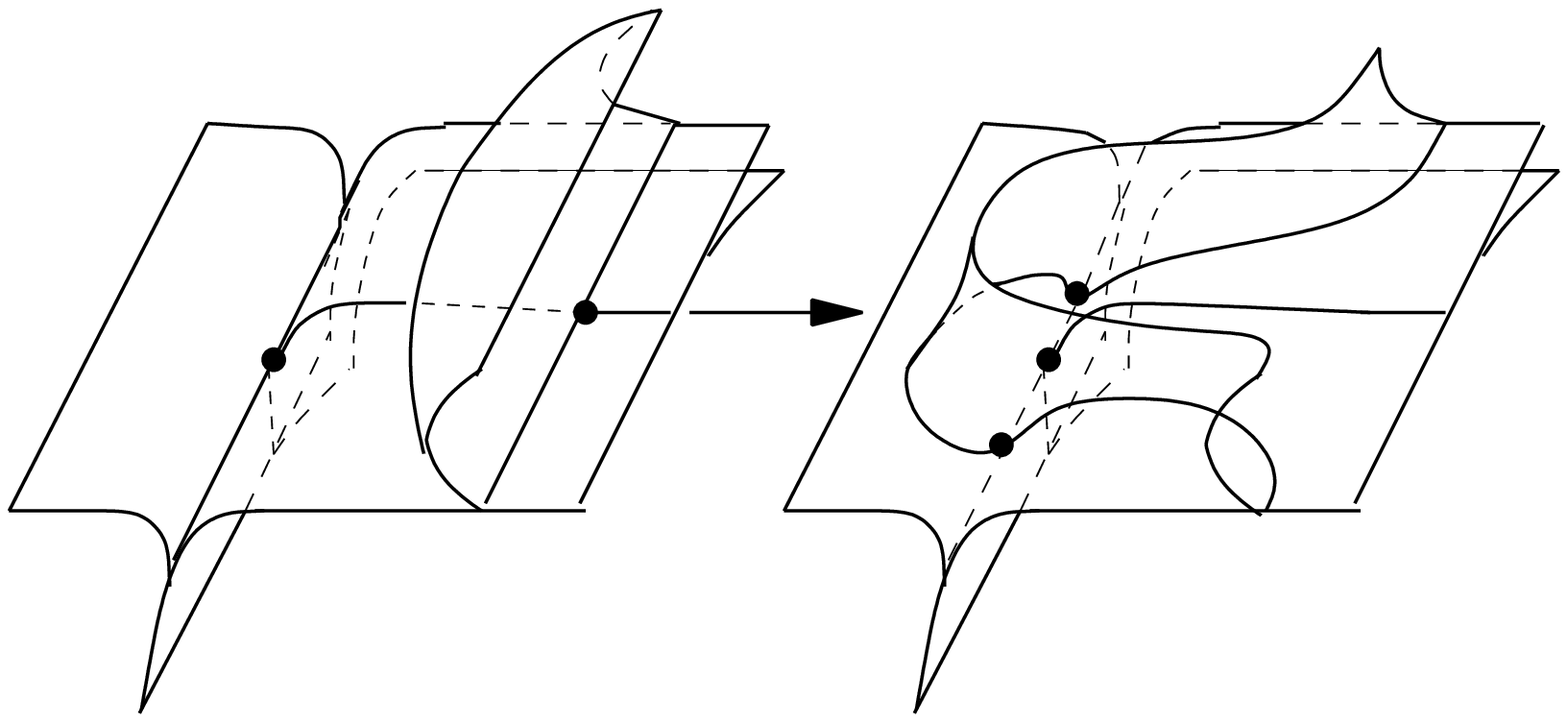}}


  \caption{In the upper part of the figure we show the $5$
 branched versions of the $2\to 3$-move called
``sliding''-moves. In the bottom part of the figure we show the
 version called ``bumping''-move. The arrow on the left indicates the
 vertical direction.
 }\label{branched MP}
\end{figure}

A branching is a kind of loss of symmetry on a
polyhedron and this is reflected by the fact that each move has
many different branched versions.
To enumerate all the possible embedded branched versions of the moves, one has
to fix any possible
orientation on the regions of the left part of Figure
\ref{2mosse} and then complete these
orientations in the right part of the figure by fixing one orientation 
on the region created by the move; by Definition \ref{branchedmove},
one obtains a branched version of a basic move when the branching
condition is satisfied both in the left and in the right part of
the figure.
Fortunately, many of the possible combinations are equivalent
up to symmetries of the pictures,
so that all the resulting local moves have been classified in \cite{BP}.
We show them in Figure \ref{branched lune} for
the lune-move and in Figure \ref{branched MP} for the $2\to 3$-move.
In these
figures we split these branched versions in two types namely the
{\it sliding}-moves and the {\it bumping}-moves; this
differentiation will be used when stating Theorem \ref{teo:calcolo
raffinato}.

\begin{figure} [h!]
   \centerline{\includegraphics[width=0.1cm]{verticaldir.eps}\includegraphics[width=11.4cm]{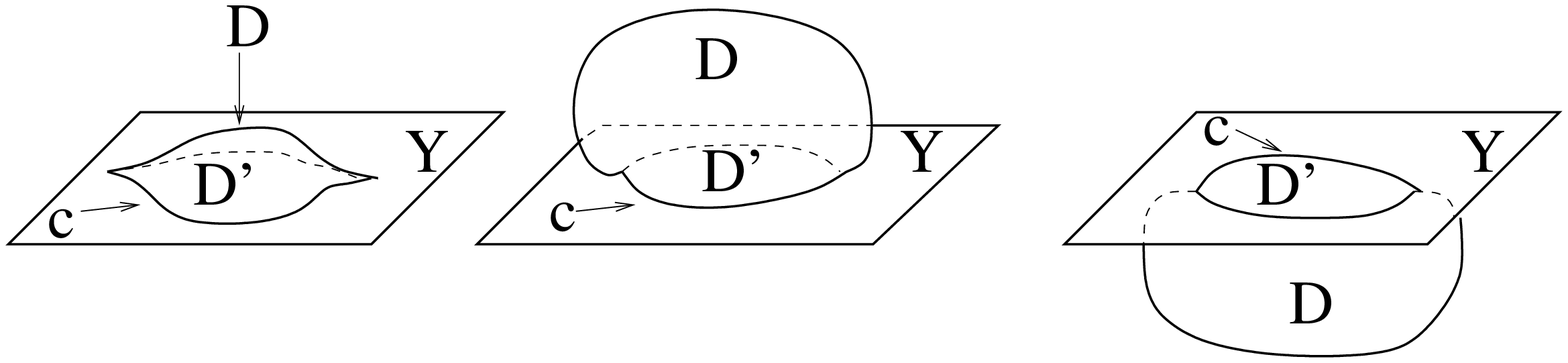}}


  \caption{In this figure we draw the $3$-versions of the embedded branched bubble-move: the left-most one is called
  ``sliding move'' and the other two ``bumping
 moves''. The arrow on the left indicates the vertical direction.
 }\label{branched bubble}
\end{figure}

We note here that in some cases 
both orientations on the region created by a move satisfy the 
branching condition, and this produces different branched versions of the same move.
This is the case of the bumping moves: for instance in the lower part of 
Figure \ref{branched MP} the orientation fixed on the disc created by
the move is such that the left most region and the disc induce 
opposite orientations on the edge separing them, so that the disc appears as a smooth 
continuation of the leftmost region out of that edge. The same phenomenon happens for leftmost sliding-moves in Figures
\ref{branched lune} and \ref{branched MP}. 
In these cases, exploiting the asymmetry of 
 the branching after the move, we say that the upper-right region ``slides over'' the
 upper-left region and creates the disc.

More in general, each instance of branched move can be viewed as a
local sliding of a region over some neighboring one: in Figures
\ref{branched lune} and \ref{branched MP} the sliding region is always
the one which, after the application of the move, is the uppermost
with respect to the depicted vertical direction. We warn the reader
that, in the following section, with an abuse of notation, we will
often call the small disc created by a positive branched-move with the
same name as the region over which the sliding has been performed;
anyhow, this abuse will always be explicitly pointed out.

It has already been proved in \cite{BP}[Chapter 3] that any lune
and $2\to 3$-move is branchable, but that there are some inverse 
lune-moves and $2\rightarrow 3$-moves which are not branchable.
Let us now analyze the branched versions of the bubble-move, referring to 
Figure \ref{fig:bubble} for the notation.

If one assigns an orientation to $D$ such that the
orientation it induces on $c$ is the same as the one induced by $D'$ on $c$
, then the bubble-move is
said to be a {\it sliding} bubble-move. 
The other case is called {\it bumping} bubble-move; as shown in Figure
\ref{branched bubble}, while performing 
an embedded bumping bubble move inside an oriented $3$-manifold,
two subversions of the move can be distinguished
depending on whether $D$ 
lies ``over'' or ``under'' $D'$; indeed the orientation of $D'$ and
the orientation of the ambient manifold allow us to distinguish an
upper and a lower face of $D'$. 

In \cite{BP}, Benedetti and Petronio proved that every orientable
$3$-manifold admits a branched spine. Moreover, in \cite{BP} and
\cite{BPco} they identified a refined topological structure which
is naturally encoded by a branching on a spine. Indeed, the
vertical vector field in Figure \ref{branching} determines a
well-defined homotopy class of vector fields on the ambient
manifold which are positively transverse to the spine, whose
orbits are properly embedded arcs and which are transverse to the
boundary of the manifold except in a finite set of simple closed
curves where they are tangent in a concave fashion. These
particular kinds of vector fields where called ``Concave
Transversing Fields'' by Benedetti and Petronio who, in \cite{BP}
and \cite{BPco} among a series of other results, proved a calculus
for manifolds equipped with these fields, which we very roughly
summarize as follows:
\begin{teo}\label{teo:calcolo raffinato}
To each embedded branched spine one can naturally associate an
homotopy class of Concave Transversing Fields; moreover, any two
branched spines encoding the same class, are connected by a
sequence of embedded branched moves of the types which in Figures
\ref{branched lune}, \ref{branched MP} and \ref{branched bubble}
are called ``sliding''-moves. If the ambient manifold is closed
then the sequence can be chosen to contain no bubble-move.
\end{teo}
The above result represents a highly refined calculus for branched
spines and was used in \cite{BPto} as a fundamental step to
produce topological invariants of homotopy classes of Concave
Transversing Fields.

\section{A calculus for branched spines of $3$-manifolds}
If one is interested in representing $3$-manifolds by means of
branched spines but is not interested in the particularly refined
structure that the branching encodes, then one needs to find a
calculus for branched spines allowing one to apply moves which
change the homotopy class of Concave Transversing Fields
represented by the spines to include all of them. We then prove the following:

\begin{teo}\label{thm:branchedspines}
 Let $M$ (with $\partial M$ possibly empty) be an oriented and compact $3$-manifold and let $P_{1}$ and
 $P_{2}$ be two branched standard spines of $M$. There exists a
 sequence of branched moves including $2\rightarrow 3$-moves, lune-moves
 and bubble-moves connecting $P_{1}$ and $P_{2}$. Moreover, the
 sequence can be constructed so that at any step the spines
 involved are standard.
\end{teo}
\begin{prf}{1}{
    By the calculus of Matveev-Piergallini, any two standard spines
    of $M$ are connected by a sequence
    of $2\rightarrow 3$-moves.
    Moreover, by a result of Y. Makovetsky (see \cite{Mak}),
    it is possible to choose two sequences of positive
    $2\rightarrow 3$-moves connecting respectively $P_{1}$ and
    $P_{2}$ to the same standard spine $P$ of $M$.
    Since these sequences are composed by positive moves, they are
    branchable and hence they connect $P_{1}$ and $P_{2}$ to two
    different branched versions of the same spine $P$. Let us call
    these two branchings on $P$ respectively $b_{1}$ and $b_{2}$.
    In what follows, we will show that, using also branched versions
    of bubble-moves, lune and $2\to 3$-moves, it is possible 
    to connect the branched spine $(P,b_1)$
    to $(P,b_2)$.

    Let $R^+_{1},.., R^+_{n}$ be the regions of $P$ oriented according to $b_1$; 
    since $M$ is
    oriented, it makes sense to speak of the upper face and of the
    lower face of a region of $P$ with respect to a branching.
    Apply a bumping bubble-move as shown in the central part of
    Figure \ref{branched bubble} to each region
    $R^+_{i}$ over which $b_1$ and $b_2$ differ, 
    so that the bubbles are attached
    along the upper face of the region with respect to $b_{1}$ 
    and call $R^-_{i}$ the new discs
    attached by the bubble-moves ($D$ in Figure \ref{branched bubble}). 
    With an abuse of notation, 
    exploiting the asymmetry of the
    bumping bubble-moves, we will call $R^+_i$ also the small disc ($D'$ in Figure 
    \ref{branched bubble}) 
    created inside $R^+_i$ by the gluing of $R^-_i$. 
    Note that the polyhedron $P'_1$ obtained after the application of these 
    moves contains
    $(P,b_1)$ as a branched sub-polyhedron and in particular for each edge or vertex of 
    $P$ there is a corresponding one in $P'_1$; moreover, $P'_1$ is 
    necessarily non standard: 
    we will sketch how to restore the standard setting in
    the end of the section. Analogously, let us call $P'_2$ the
    polyhedron obtained by applying the above procedure using the
    branching $b_2$ to choose the upper faces of the regions. 
    It is clear that $P'_1$ and $P'_2$
    are branched version of the same polyhedron $P'$; we
    keep calling $b_1$ and $b_2$ their branchings.

    The idea of the proof is to slide the regions $R^-_i$ in $P'_1$ to a 
    final position where the roles of
    $R^+_{i}$ and of $R^-_{i}$ are exchanged and so, in particular, each
    $R^+_{i}$ appears as a bubble applied to the center of $R^-_{i}$, on
    its upper face: this connects $P'_1$
    and $P'_2$ by a sequence of branched moves (not including bubble-moves).
    
    Here and in what follows we will use the
    natural identification induced by a branched move between the
    regions of a polyhedron before and after the move as explained in
    the preceding section. Moreover, for the sake of simplicity, we suppose that 
    the edges and the vertices we deal with are touched by distinct regions; 
    this is not true in general but the proof can be easily adapted to the general case using 
    local names for the regions around the edges and the vertices instead of the names 
    $R^\pm_*$.    
    We split the proof in two
    steps:
    \begin{enumerate}
    \item For each edge $e$ of $P$ we use positive
    lune-moves near the corresponding edge $f$ in $P'_1$ to create an edge such that if $e$ is
    touched by $R^+_i$, $R^+_j$ and $R^+_k$ in $P$ then the new edge is touched by
    $F(R^+_i)$, $F(R^+_j)$ and $F(R^+_k)$, where $F(R^+_*)=R^+_*$ if $b_1$ equals $b_2$ on 
    $R^+_*$ in $P$ and is $R^-_*$ otherwise. This
    sequence of moves modifies $P'_1$ creating some additional
    singular locus and hence producing a polyhedron no longer
    homeomorphic to $P'$, but since it acts only near the center
    of the edges of $P'_1$, the vertices of $P'_1$ and hence of $P$ are
    identified with a subset of the vertices of the resulting polyhedron.
    \item For each vertex $v$ of $P$, we apply a sequence of branched
    moves to the corresponding vertex $w$ in the polyhedron obtained after Step $1$ to
    create a vertex 
    such that if $R^+_i$, $R^+_j$, $R^+_k$, $R^+_l$, $R^+_m$ and $R^+_n$ are
    the regions of $P'_1$ touching $v$ then the vertex is touched by the
    regions $F(R^+_i)$, $F(R^+_j)$, $F(R^+_k)$, $F(R^+_l)$, $F(R^+_m)$ and 
    $F(R^+_n)$ (where,
    again, $F(R^+_*)=R^+_*$  if $b_1$ and $b_2$ are equal in $P$ on $R^+_*$ and is 
    $R^-_*$ 
    otherwise).
    The sequence of moves we apply has also the effect 
    of eliminating the extra singular
    locus created during the first step so that the final polyhedron
    we get is $P'_2$.

    \end{enumerate}

    {\bf Step 1.}
    First of all, note that there are exactly $6$ possible different
    branchings on the neighborhood of an edge of a spine 
    (of the $8$ possible $3$-uples of orientations on the regions
    touching the edge, two are to be excluded
    since of the branching condition). 
    Moreover, each branching on an edge
    together with the orientation of $M$,
    produces a cyclic ordering of the regions touching it
    and a notion of up and down near it, so that we can define as shown in
    Figure \ref{fig:types} the regions of type $1$, $2$ and $3$ with
    respect to the branching on the edge. 

    \begin{figure}[htbp]
    \centerline{\includegraphics[width=0.1cm]{verticaldir.eps}\includegraphics[width=8.4cm]{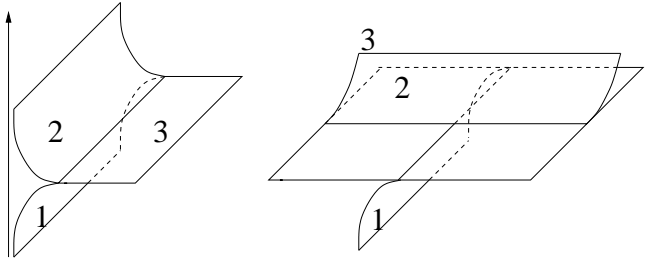}}

    \caption{In this picture we show how to determine the regions of
    type $i=1,2,3$ near an edge and near a vertex. In the picture  we
    use the right-hand rule and vertical direction drawn on the left
    to determine the orientations of the regions. }
    \label{fig:types}
    \end{figure}
    Let $e$ be an edge of $P$ where $b_{1}$ and $b_{2}$ are
    different, $f$ the 
    corresponding edge in $P'_1$ and let respectively $R^+_i$, $R^+_j$ and $R^+_k$ be the
    regions of type $1$, $2$ and $3$ with respect to the branching
    $b_1$ on $e$. 
    
    Suppose that $R^+_i$ is the only region touching $e$ on which $b_1$ and $b_2$ 
    differ;
    we will now exhibit a sequence of $2$ positive bumping
    lune-moves near $f$ in $P'_1$ which creates an edge 
    touched by $R^-_i$, $R^+_j$ and
    $R^+_k$ and some additional singular locus which will be eliminated
    in Step $2$.  
    
    Since $R^+_i$ is oriented in the opposite way by $b_1$
    and $b_{2}$, to get $P'_1$ from $P_1$
    we applied a bumping bubble-move over this region and created a
    disc $R^-_{i}$ on the upper (with respect to $b_{1}$) face of it.
    Apply a bumping lune-move to slide $R^-_{i}$ over a
    little disc contained in $R^+_j$ 
    as shown in the
    upper part of Figure \ref{fig:edge1}. This splits $f$
    adding two new vertices and creates a new small
    disc and a new singular edge $f'$ touched by $R^-_i$, $R^+_j$ and the disc
    itself. Apply now another bumping lune-move to slide $R^+_i$ over
    $R^-_i$ passing through the disc as shown in the lower-left part of
    Figure \ref{fig:edge1}. The edge $f'$ is split by two
    new vertices and a new small disc is created. 
    With an abuse of notation we will call this
    new small disc $R^-_i$; indeed it can be viewed as a smooth
    continuation of $R^-_i$ over the edge separing them. 
    After this move, the small straight edge connecting the two vertices 
    in the lower-right part of Figure \ref{fig:edge1} is touched
    by the regions $R^-_i$, $R^+_j$ and $R^+_k$: we
    obtained the singular edge we were searching for. 
    The side-effect of the sequence of moves just described is to create
    two pairs of new vertices positioned symmetrically with respect to
    the center of the edge: we will show in Step $2$ how to eliminate these vertices.

    \begin{figure}[htbp]
    \centerline{\includegraphics[width=8.4cm]{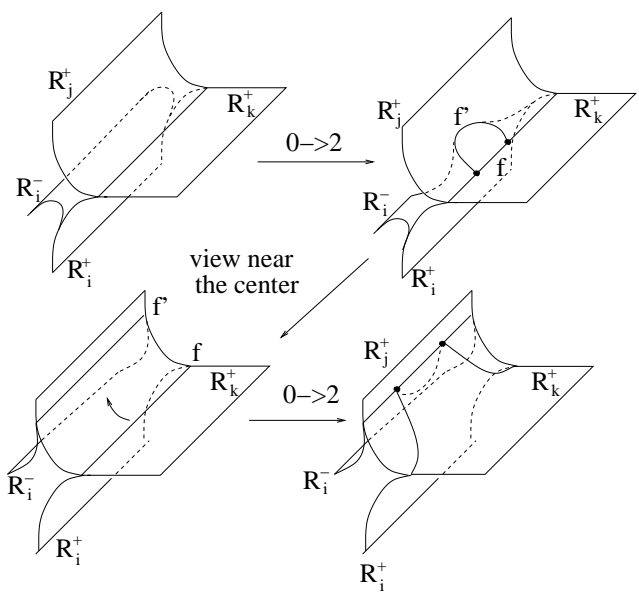}}
    \caption{In this picture we show the sequence of moves uses to exchange near the 
    center of an
    edge of $P'_1$ a 
    region $R^+_i$ of type $1$ with respect to
    $b_1$ with its companion $R^-_i$ . 
    The first move is a bumping lune-move shifting $R^-_{i}$ over
    $R^+_{i}$; then we concentrate on the central part of the picture and
    we apply another bumping lune-move to slide $R^+_{i}$ over $R^-_{i}$.}
    \label{fig:edge1}
    \end{figure}

    \begin{figure}[htbp]
    \centerline{\includegraphics[width=8.4cm]{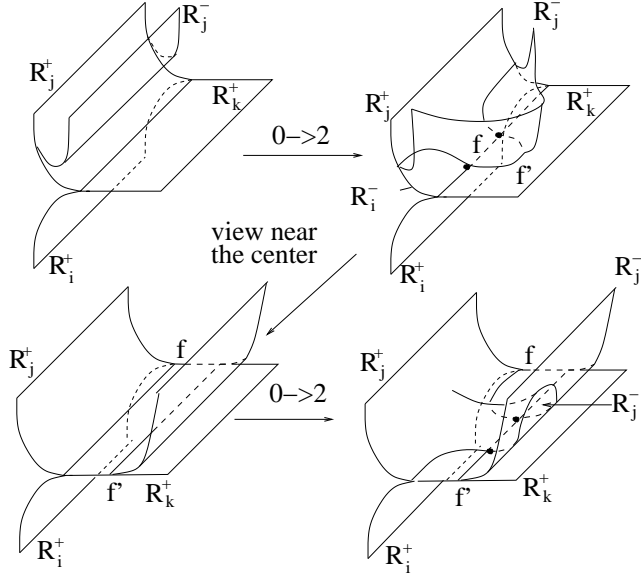}}
    \caption{In this picture we show how to revert the orientation of a
    region of type $2$ near the center of an edge. In the first step we
    apply a sliding lune-move to push $R^-_{i}$ over $R^+_{i}$, then we
    concentrate on the center part of the picture and we apply another
    sliding lune-move to push $R^+_{i}$ over $R^-_{i}$.}
    \label{fig:edge2}
    \end{figure}

    Suppose now that $b_1$ and $b_2$ differ only on $R_j$.
    We will now exhibit a sequence of $2$ sliding lune-moves 
    near $f$ in $P'_1$
    creating a new edge $f'$ touched by $R^+_i$, $R^-_j$ and
    $R^+_k$.  Apply a sliding lune-move to push $R^-_{j}$ onto a
    little disc contained in $R^+_k$ and containing the center of $f$ in its boundary.
    The boundary of the
    little disc is formed by $f$ and another small singular edge $f'$
    (see the upper part of Figure \ref{fig:edge2}).
    Then, we apply another sliding lune-move to push $R^+_{j}$
    over $R^-_{j}$ near the center $f'$ passing over the small
    disc. This creates another small disc ``contained'' in $R^-_j$
    which, with an abuse of
    notation, we will call $R^-_j$. The sequence creates the edge
    touched by $R^+_i$, $R^-_j$ and $R^+_k$ we were searching for.
    Again, as a side-effect of the sequence, two pairs of new vertices disposed
    symmetrically with respect to the edge appear.
    
    Let us now show how to conclude Step $1$ in the general case around $e$. 
    Starting from $f$, we want to produce an edge 
    in $P'_1$ touched by the regions
    $F(R^+_*)$ (where $F(R^+_*)=R^+_*$ iff $b_1$ and $b_2$ do not
    differ on $R^+_*$). To do this, we ``substitute" the regions one by one using the
    above sequences of moves.  
    Consider the following algorithm: 
    \begin{enumerate}
      \item Let $X$ be the region of type $1$ with respect to the
      branching on $f$. If $X=F(R^+_*)$ for some $R^+_*$ in $P'_1$,  go to $2$, otherwise 
      apply the sequence of Figure
      \ref{fig:edge1} to $f$: this produces a new branched edge 
      which we keep calling $f$ whose regions
      are still indicized by the three-uple $\{i,j,k\}$ and
      in which one more region is of the form $R^-_*=F(R^+_*)$. Go to
      $2$. 
      \item Let $Y$ be the region of type $2$ with respect to the
      branching on $f$. If $Y=F(R^+_*)$ for some $R^+_*$ in $P'_1$ stop, otherwise 
      apply the sequence of Figure
      \ref{fig:edge2} to $f$: this produces a new branched edge 
      which we keep calling $f$ whose regions
      are still indicized by the three-uple $\{i,j,k\}$ and
      in which one more region is of the form $R^-_*$; go back to $1$.   
    \end{enumerate}
    Since the regions around $f$ are indicized in $\{i,j,k\}$, we can ``pull-back" on $e$ the 
    branching of $f$ and call it $b(f)$. 
    To compare $b(f)$ and $b_1$ on $e$ we stipulate that they differ on a region $R^+_x$ 
    touching $e$ iff the corresponding region touching $f$ is $R^-_x$. Consequently, 
    $b(f)$ and $b_2$ 	 
    differ on a region $R_x$ near $e$ iff the corresponding region touching $f$ is $R^+_x$ 
    and  $b_1$ and $b_2$ differ on $R_x$. 
    
    We claim that when the above algorithm stops $b(f)=b_2$ on $e$. Indeed they cannot differ 
    on a region of type $1$ or $2$ near $e$ with respect to $b(w)$ since otherwise the 
    algorithm would not have stopped.  But then, they cannot differ on the remaining region since 
    otherwise $b_2$ would not satisfy the branching condition on $e$ (by construction 
    $b(w)$ does).
    
    Roughly speaking, by applying the above procedure to $P'_1$ for each edge of
    $P$, one obtains
    a new branched spine which, near the centers of some edges 
    is branched according to $b_{2}$, near the vertices corresponding
    to those of $P$ is still branched
    according to $b_{1}$ and which contains some additional singular locus 
    which will delete in the next Step.

    {\bf Step 2.}
    Let $w$ be a vertex of the branched polyhedron obtained after Step $1$ and
    corresponding to a vertex $v$ of $P$. 
    A branching canonically equips each edge of a spine with an orientation (the one 
    induced by two of the $3$ regions touching the edge) and this allows us to identify $3$ 
    particular regions around each branched vertex. These regions are those 
    which, in the neighborhood of the vertex, orient both the edges in their boundary 
    positively,  and are shown in Figure \ref{fig:types}. We call 
    them regions of type $1$, $2$ and $3$; we tell
    regions of type $1$ from those of type $3$ by means of the orientation of $M$:
    the type $1$ is the lower one with respect to the positively
    oriented normal to the spine in the vertex. 
    
    Suppose by now that $b_1$ and $b_2$ differ near $v$ only on the region $R^+_i$ of 
    type $1$
    w.r.t $b_1$ and let us exhibit the sequence concluding Step $2$ on
    $v$ in this case. 
   
    Let $e_{1}$ and $e_{2}$ be the two edges of $P$ touching $v$ and
    contained in the boundary of $R^+_{i}$. Note that $R^+_{i}$ is of
    type $1$ w.r.t $b_1$ also for $e_1$ and $e_2$ (see Step $1$) and
    hence during Step $1$ we applied the 
    sequence of moves of Figure \ref{fig:edge1} to the edges $f_1$ and $f_2$ in $P'_1$ 
    corresponding to $e_1$ and $e_2$. 
    In particular, the first
    moves applied on these two edges during Step $1$ are bumping
    lune-moves shifting $R^-_{i}$ over small discs contained in the regions of type
    $2$ respectively for $f_{1}$ and for $f_{2}$ and which, near $w$
    appear as horizontal. Then near $w$ we see the pattern
    described in the upper-left part of Figure \ref{fig:type1}.
    Apply an inverse bumping $2\rightarrow 3$-move
     and complete
    the shifting of $R^-_{i}$ near $w$ (see the upper-right part of Figure
    \ref{fig:type1}). Then, applying a
    bumping $2\rightarrow 3$-move we slide $R^+_i$ over a small disc on
    $R^-_i$ which, with an abuse of notation, we will call $R^-_i$. 
    Finally, since the last moves of the sequence of Figure
    \ref{fig:edge1} applied in Step $1$ 
    are two bumping lune-moves near $f_1$
    and $f_2$ sliding $R^+_i$ over $R^-_i$, we are in the situation
    depicted in the lower-left part of the same figure. We can then
    apply two inverse bumping lune-moves to eliminate the two small
    discs horizontal in the lower-left part of the figure and slide
    the region $R^+_{i}$ over
    $R^-_{i}$ as shown in the lower-right part of the figure. 
    That way, the vertices which had been created during Step $1$ on
    the edges $f_1$ and $f_2$ disappear and we created a new vertex
    touched by the same regions as $v$ with the only exception of $R^+_i$
    which has been substituted by $R^-_i$. This concludes Step $2$ in
    this case.

    \begin{figure}[htbp]
    \centerline{\includegraphics[width=10.4cm]{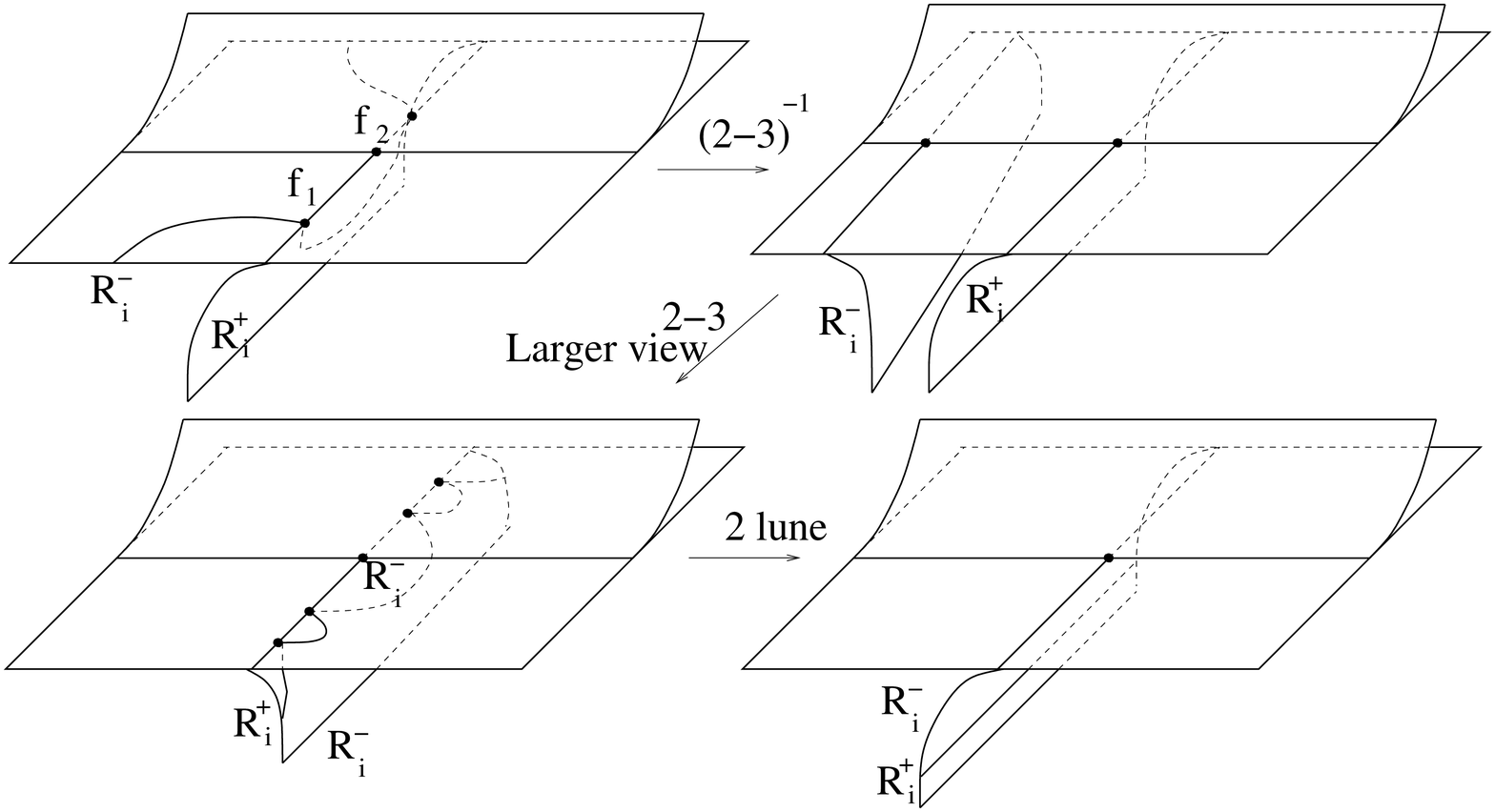}}
    \caption{In this picture we show the sequence used in the case
      when the branchings on a vertex differ on a region of type $1$
    with respect to the first branching.}
    \label{fig:type1}
    \end{figure}
 Let us now suppose that $b_1$ and $b_2$ differ near $v$ only in
    the region of type $2$ w.r.t $b_1$ and let $R^+_j$ be this region.

    Let again $e_{1}$ and $e_{2}$ be the two edges of $P$ in the boundary of
    $R^+_{j}$ touching $v$ and $f_1$ and $f_2$ the corresponding edges in $P'_1$ .
    Note that, w.r.t $b_1$, $R^+_{j}$ is of type $1$ for one of
    the two edges (say $e_{1}$) and of type $2$ for the other one
    ($e_{2}$).
    Moreover, note that the region near $v$ which is of type $1$ for
    $e_{2}$ is also of type $1$ for $v$, hence, by hypothesis,
    this region in $P$ is oriented
    the same way by $b_{1}$ and $b_{2}$.
    Then the moves applied on $f_{1}$ during Step $1$ are those of the
    sequence of Figure \ref{fig:edge1}: in particular the first one is 
    a bumping lune-move 
    sliding the region $R^-_{j}$ over the center of the edge. On $f_{2}$,
    the moves of Step $1$ are those of the sequence of Figure
    \ref{fig:edge2} and the first one is a sliding lune-move pushing
    $R^-_{j}$ near the center of the edge over a small disc 
    contained in the region of type $3$ for $f_{2}$.
    Hence, near $w$, the polyhedron after Step $1$ appears as shown in
    the upper-left square of Figure \ref{fig:type2}.
    \begin{figure}[htbp]
    \centerline{\includegraphics[width=10.4cm]{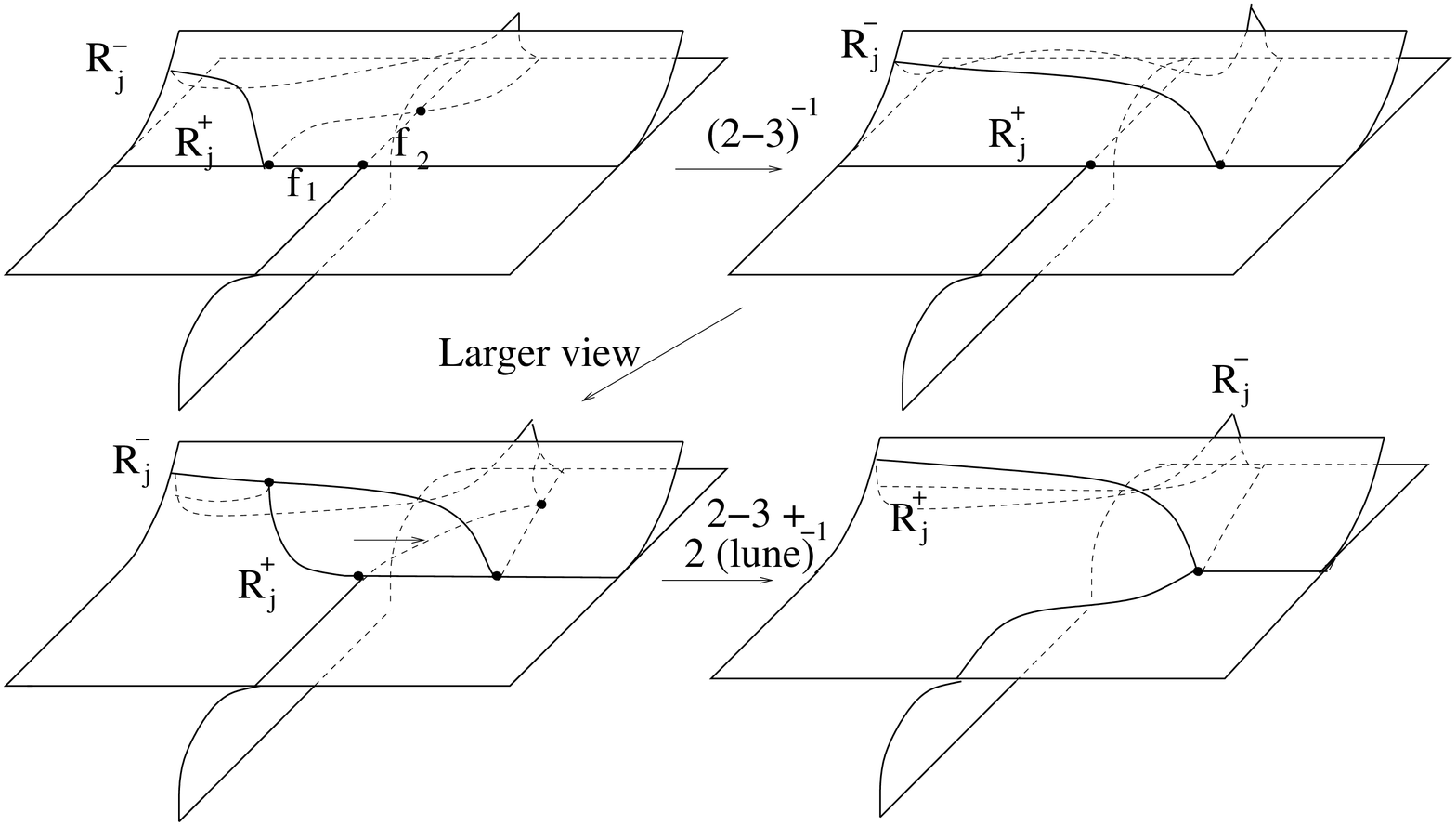}}
    \caption{In this picture we show the sequence used in the case
      when the branchings on a vertex differ on a region of type $2$
    with respect to the first branching. 
    Note that, for the sake of clarity, in the above picture we do not the draw some regions
    and we limit ourselves to outline their boundary curves.}
    \label{fig:type2}
    \end{figure}

    By applying an inverse sliding $2\rightarrow
    3$-move, one shifts $R^-_{j}$ over $w$ as shown in the upper-left
    part of the figure.  
    Finally, since the last moves of the sequences applied during Step $1$ 
    on the edges $f_1$ and $f_2$ slid $R^+_j$
    over $R^-_j$ we are left in the case shown in the lower-left part
    of the figure. Then we apply a $2\rightarrow 3$-move
    to slide $R^+_j$ over a small disc contained in $R^-_j$ 
    as indicated by the small arrow in the lower-left part of the figure; with an
    abuse of notation we call this disc $R^-_j$.  Then we perform
    two inverse lune-moves to eliminate the two small discs
    visible in the lower-left part of the figure and complete the
    sliding of $R^+_{j}$ over $R^-_{j}$. This completes Step $2$ in this case. 

    Until now we showed how to conclude Step $2$ when $b_1$ and
    $b_2$ differ on a region of type $1$ or $2$ near $v$, 
    Now we examine
    the case when $b_1$ and $b_2$ differ only on the region
    $R^+_{k}$ of type
    $3$ w.r.t $b_1$ near $v$.

    Let again $e_{1}$ and $e_{2}$ be the two edges in
    the boundary of $R^+_{k}$ and $f_1$ and $f_2$ the corresponding edges in $P'_1$,
    so that $e_{2}$ is also in the boundary of
    the region which is of type $2$ with respect to $b_1$ on $v$. Note that
    $R^+_{k}$ is of type $2$ for both the edges w.r.t. $b_1$. Moreover, by hypothesis,
    the regions of type $1$ and $2$ w.r.t. $b_1$  near
    $v$ are oriented the same way by $b_{1}$ and $b_{2}$,
    this implies that also the region near $v$ which contains in its
    boundary $e_{1}$ and orients it in the same way as $R^+_{k}$ is
    oriented the same way by $b_{1}$ and $b_{2}$ (otherwise $b_2$
    would not define a branching). Hence, during Step $1$ we applied
    both on $f_1$ and on $f_2$ the sequence of moves of Figure \ref{fig:edge2}.
    Then, near $w$, we are in the situation depicted in the upper-left part of
    Figure \ref{fig:type3}. We
    apply an inverse sliding $2\rightarrow 3$-move near $w$ to slide $R^-_{k}$
    over the horizontal plane (see the upper-right part of the
    figure). 
    Then we apply a sliding $2\rightarrow 3$-move to slide $R^+_k$ over
    a small disc contained $R^-_k$ near the vertex. With an abuse of
    notation we will call this disc $R^-_k$. Since in the sequence of
    moves applied during Step $1$ over each edge a
    sliding lune-move has been performed to slide $R^+_k$ over $R^-_k$,
    we are now in 
    the situation depicted in the lower-left square of the figure.
    To finish, we then apply 
    two inverse sliding lune-moves to slide $R^+_{k}$ over $R^-_{k}$ and
    eliminate the two small discs present in the lower left part of
    the figure.
 \begin{figure}[htbp]
    \centerline{\includegraphics[width=10.4cm]{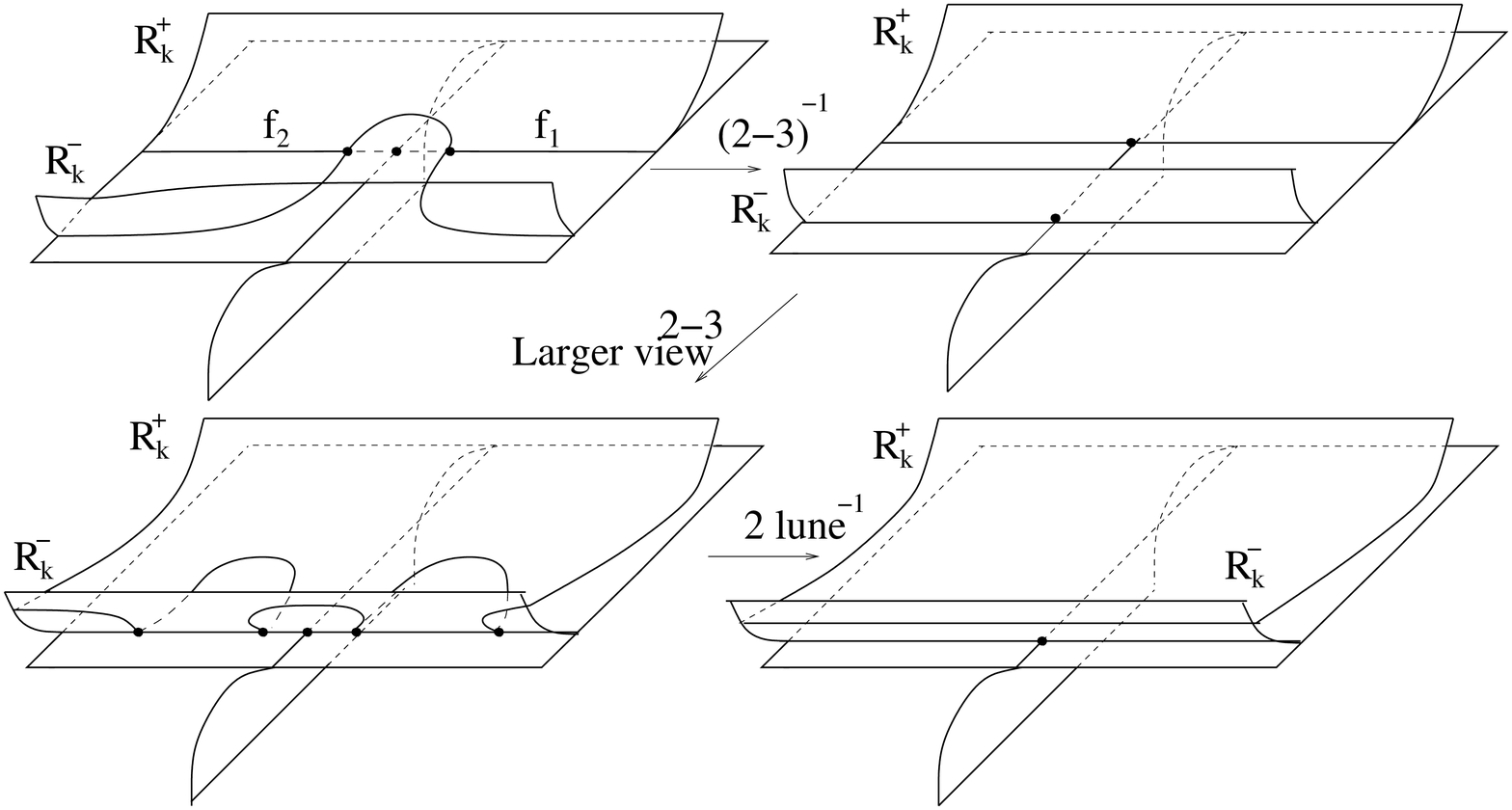}}
    \caption{In this picture we show the sequence used in the case
      when the branchings on a vertex differ on a region of type $3$. }
    \label{fig:type3}
    \end{figure} 
    
    We now claim that the sequences shown in Figures \ref{fig:type1}, 
    \ref{fig:type2} and \ref{fig:type3} are sufficient to complete
    Step $2$ in the general case when $b_1$ and $b_2$ differ on more
    than one region near $v$. 
 
    Now let $R^+_i$,$R^+_j$,$R^+_k$,$R^+_l$,$R^+_m$ and $R^+_n$ be the regions of 
    $P$
    around $v$. In the vertex $w$ of 
    the polyhedron obtained
    after Step $1$ and corresponding to $v$, we want to ``substitute''  each region
    $R^+_*$ with $F(R^+_*)$ (where $F(R^+_*)=R^+_*$ iff $b_1$ and $b_2$ do not
    differ on $R^+_*$), that is apply suitable sequence
    of moves which produce a new vertex touched by the six regions
    $F(R^+_*)$. To do this we substitute the regions one by one using the
    sequences exhibited above.  
    Consider the following algorithm: 
    \begin{enumerate}
      \item Let $X$ be the region of type $1$ with respect to the
      branching on $w$. If $X=F(R^+_*)$ go to $2$, otherwise 
      apply the sequence of Figure
      \ref{fig:type1} to $w$: this produces a new branched vertex 
      which we keep calling $w$ whose regions
      are still indicized by the six-uple $\{i,j,k,l,m,n\}$ and
      in which one more region is of the form $R^-_*$. Go to
      $2$. 
      \item Let $Y$ be the region of type $2$ with respect to the
      branching on $w$. If $Y=F(R^+_*)$ go to $3$ otherwise 
      apply the sequence of Figure
      \ref{fig:type2} to $w$: this produces a new branched vertex 
      which we keep calling $w$ whose regions
      are still indicized by the six-uple $\{i,j,k,l,m,n\}$ and
      in which one more region is of the form $R^-_*$; go back to $1$. 
      \item Let $Z$ be the region of type $1$ with respect to the
      branching on $w$. If $Z=F(R^+_*)$ stop otherwise apply the sequence of Figure
      \ref{fig:type3} to $w$: this produces a new branched vertex 
      which we keep calling $w$ whose regions
      are still indicized by the six-uple $\{i,j,k,l,m,n\}$ and
      in which one more region is of the form $R^-_*$. End. 
    \end{enumerate}
    Since $w$ is branched and the regions touching it are indexed
    in $\{i,j,k,l,m,n\}$, we can ``pull back" on $v$ the branching on $w$ which we will call 
    $b(w)$ and compare it 
    with $b_1$ and $b_2$:  $b(w)$ differs from $b_1$ on a 
    region $R^+_x$ if and only if $w$ is touched by $R^-_x$.  Analogously $b(w)$ and $b_2$ 
    differ on a region $R^+_x$ iff $b_1$ and $b_2$ differ on it and $w$ is touched by the 
    region $R^+_x$. 
    
    We claim that, when the algorithm stops, the regions of the form $R^-_*$ touching $w$
    are exactly those over which $b_1$ and $b_2$ differ and then $b_2=b(w)$ 	
    near $v$. 
    Indeed, the above algorithm stops when the
    regions of type $1$, $2$ and $3$ around $w$ are all of the form
    $F(R^-_*)$; clearly, this is 
    achieved in at most $6$ steps. 
    Then, in particular, $b_2$ 
    and $b(w)$ define branchings coinciding on the regions of type $1$,$2$ and $3$ w.r.t. 
    $b(w)$. This implies that $b_2$ and $b(w)$ are equal on all the regions since if $b_2$ 
    differed from $b(w)$ on another region then the branching condition would not be 
    satisfied by $b_2$.
    It is important to note that the algorithm above can be followed until its end since 
    during Step $1$ we used an algorithm based on the edges of $P$ producing exactly 
    the extra vertices eliminated by the sequences of Figures \ref{fig:type1}, \ref{fig:type2} 
    and \ref{fig:type3} and ordered compatibly. This concludes Step $2$.
    
    Now let us note that after Step $2$, the polyhedron one obtains is homeomorphic to 
    $P'$ (during Step $2$ we eliminated all the extra structures constructed during Step $1$), and its 
    branching coincides with the one given by $b_2$ since, 
    using embedded branched moves, we exchanged the roles of the 
    regions $R^+_*$ with those of the regions $R^-_*$ having the opposite orientation.  
    We are now left to prove that the sequence of moves we used can be improved to a 
    sequence passing only through standard polyhedra: we limit ourselves to 
    sketch the idea of the proof. For each region $R^+_i$ in $P$ let us choose an edge 
    $e_i\subset Sing(P) $ touched by $R^+_i$. Then at the beginning of the proof, if 
    $b_1(R^+_i)\neq b_2(R^+_i)$, instead of applying a bubble move inside $R^+_i$ we 
    apply a standard bubble-move (whose non branched version is shown in the right part 
    of Figure \ref{fig:bubble}) gluing $R^-_i$ along $e_i$ slightly aside from the center of 
    the edge and on the upper side of $R^+_i$ (with respect to $b_1$). The curve 
    $\partial R^-_i$ will bound a disc in $P$ formed by the union of two discs, one 
    ``contained" in $R^+_i$ and the other, call it $D_i$ coming from another region 
    touching $e_i$.  Now apply Step $1$ and then
    apply an inverse lune-move to slide $R^-_i$ out of $D_i$. It can be checked 
    that the polyhedron one obtains is equal to the result of Step $1$ obtained by following 
    the proof in 
    the non standard case, but never passes through non standard polyhedra so that one 
    can proceed with Step $2$. To ensure that also during Step $2$ one does not produce 
    non standard polyhedra, it is sufficient to observe that the only moment when such a 
    polyhedron can be obtained is when one concludes one of the sequences of Figures 
    \ref{fig:type1}, \ref{fig:type2} and \ref{fig:type3} through the pair of inverse lune moves 
    and the region to be slid, after the move, becomes a bubble over another region of the 
    polyhedron. In this case it is sufficient to perform just one of the two inverse lunes and 
    proceed with Step $2$ (which will no more involve that region).   
     }\end{prf}

\noindent
Scuola Normale Superiore\\
Piazza dei Cavalieri 7, 56127 Pisa, Italy\\
f.costantino@sns.it\\


\begin{thebibliography}{99}

\bibitem{BB}\textsc{Baseilhac, S. and Benedetti, R.}, ``{QHI}, 3-manifolds scissors congruence classes and the volume conjecture'',
Geom. Topol. Monogr. Invariants of knots and $3$-manifolds 4
(2001), 13-28.

\bibitem{BB2}\textsc{Baseilhac, S. and Benedetti, R.}, ``{QHI} Theory II, Dilogarithmic and Quantum Hyperbolic Invariants of 3-Manifolds with {$PSL(2,\mc)$}-Characters
'', xxx.arxiv.org/math.GT/0211061 (2002).

\bibitem{BB3}\textsc{Baseilhac, S. and Benedetti, R.}, ``Quantum Hyperbolic Invariants Of 3-Manifolds With {PSL(2,$\mc$)}-Characters
'', xxx.arxiv.org/math.GT/0306280 (2003).

\bibitem{BB4}\textsc{Baseilhac, S. and Benedetti, R.}, ``Classical And Quantum Dilogarithmic Invariants of Flat {PSL(2,$\mc$)}-Bundles Over 3-Manifolds
'', xxx.arxiv.org/math.GT/0306283 (2003).

\bibitem{BP}\textsc{R.~Benedetti, C.~Petronio}, ``Branched
Standard Spines of $3$-manifolds'', Lecture Notes in Math., vol.
1653, Springer-Verlag, Berlin-Heidelberg-New York, 1997.

\bibitem{BPto}\textsc{R.~Benedetti, C.~Petronio}, ``{Reidemeister-Turaev} torsion of $3$-dimensional {Euler} structures with simple boundary
  tangency and pseudo-legendrian knots'', Manuscripta Mathematica 106 (2002) 1, 13-61.

\bibitem{BPco}\textsc{R.~Benedetti, C.~Petronio}, ``Combed $3$-manifolds with concave boundary, framed links, and pseudo-{Legendrian} links'',
J. Knot Theory Ramifications 10 (2001) 1, 1-35.

\bibitem{DPS}\textsc{J.L.~Dupont, C-H.~Sah}, ``Scissors congruences II'',
J. Pure and App. Algebra 44 (1987), 137-164.

\bibitem{Ma}\textsc{S.~Matveev}, ``Transformations of special
spines and the Zeeman conjecture'', Math. USSR Izvestia 31 (1988),
423-434

\bibitem{Mak}\textsc{A.Yu.~Makovetski, A. Yu.}, ``Transformations of special spines and special polyhedra'',
Math. Notes 65 (1999), 295-301.

\bibitem{Pi} \textsc{R. ~Piergallini}, ``Standard moves for standard
polyhedra and spines'', Rend. Circ. Mat. Palermo 37 (1988), suppl.
18, 391-414.

\bibitem{TV} \textsc{V. ~Turaev, O. ~Viro }, ``State sum invariants of $3$-manifolds and $6j$-symbols'', Topology 31 (1992), 865-904.


\end{thebibliography}
\end{document}